\documentclass[12pt,twopage]{amsart}
\usepackage{euscript}
\usepackage{amsfonts}
\usepackage{amssymb}
\usepackage{amsmath}
\usepackage{amsthm}
\usepackage{amscd}
\usepackage{xypic}
\usepackage{enumerate}

\usepackage{amsopn}
\usepackage[english]{babel}
\usepackage[latin1]{inputenc}
\swapnumbers

\theoremstyle{plain}
\newtheorem{thm}{Theorem}[section]
\newtheorem{cor}[thm]{Corollary}

\newtheorem{prop}[thm]{Proposition}
\theoremstyle{definition}
\newtheorem{defn}{Definition}[section]
\theoremstyle{remark}
\newtheorem{rem}{Remark}[section]
\numberwithin{equation}{section}

\newcommand{\Z}{\mathbb Z}

\newcommand{\C}{\mathbb C}

\newcommand{\G}{\mathbb G}
\newcommand{\PP}{{\mathbb P}}

\newcommand{\OO}{\mathcal{O}}

\newcommand{\II}{\mathcal{I}}

\newcommand{\LL}{\mathcal{L}}
\newcommand{\VV}{\mathcal{V}}

\DeclareMathOperator{\HH}{H} \DeclareMathOperator{\hh}{h}
 
\DeclareMathOperator{\im}{Im} \DeclareMathOperator{\rk}{rk}
\DeclareMathOperator{\Ker}{Ker} 
\DeclareMathOperator{\coker}{Coker} 
 
\DeclareMathOperator{\T}{T}
\DeclareMathOperator{\s}{S}

\DeclareMathOperator{\points}{points}

\begin{document}

\title{The Horrocks-Mumford bundle restricted to planes}

\author{Ada Boralevi}

\address{Dipartimento di Matematica ``U. Dini'', Universit\`{a} degli Studi di
Firenze, Viale Morgagni 67/a, I-50134 Florence, Italy}

\email{boralevi@math.unifi.it}

\date{}



\maketitle

\begin{abstract}
We study the behavior of the Horrocks-Mumford bundle $F_{HM}$ when restricted to a plane $\PP^2 \subset \PP^4$, looking for all possible minimal
free resolutions for the restricted bundle. To each of the 6 resolutions (4 stable and 2 unstable) we find, we then associate a subvariety of
the Grassmannian $\G(2,4)$ of planes in $\PP^4$. We thus obtain a filtration of the Grassmannian, which we describe in the second part of this
work.
\end{abstract}

\section{Introduction}
The Horrocks-Mumford bundle is a stable rank-2 complex vector bundle
on $\PP^4$, and at the present state, it is -up to twist by line
bundles and finite pullbacks- the only one of its kind known to be
undecomposable (cf. Hartshorne's conjecture, \cite{Ha74}).\\
\indent It can be defined as the cohomology of the following monad:
\begin{equation}\label{monadeHM}
0\rightarrow 5{\OO_{\PP^4}(-1)} \xrightarrow{B}
2{\Omega}^2_{{\PP}^4}(2) \xrightarrow{A} 5{\OO_{\PP^4}} \rightarrow
0.
\end{equation}
Once we have equipped the 5-dimensional complex vector space $V$
with a basis $\{e_i\}_{i \in \Z_5}$, $A=(a_{ij})$ is a $2 \times 5$
matrix of 2-forms (see \cite{BHM} or \cite{DS}):
$$
A=\left(
    \begin{array}{c}
      a_{0i} := e_{i+2} \wedge e_{i+3} \\
      a_{1i} := e_{i+1} \wedge e_{i+4}\\
    \end{array}
  \right)
$$
and $B=^t(A\cdot Q)$, where $Q$ is just the matrix $Q=\left(
                                                        \begin{array}{cc}
                                                          0 & 1 \\
                                                          -1 & 0 \\
                                                        \end{array}
                                                      \right)
$.\\
The construction of the bundle is
due to Horrocks and Mumford, who discovered it in 1972 \cite{HM}.
Nevertheless via the Hartshorne-Serre correspondance its existence
is hand in glove with that of degree 10 abelian varieties in
$\PP^4$, that had already been
proved by Comessatti \cite{Co} in 1916.\\
\indent This work aims at studying the behavior of the HM-bundle $F_{HM}$
when restricted to a plane $\PP^2 \subset \PP^4$.\\
Stable and unstable planes (the latter are also called jumping) are already known by \cite{BHM}. Here we
want to study minimal free resolutions for the restricted bundle.\\
We find out that there are 6 possible minimal resolutions
$(\underline{a},\underline{b})$, of the form:
\begin{equation}\label{risoluzione}
 0\rightarrow \oplus_{i=1}^{k} \OO_{\PP^2}(-a_i)\xrightarrow{\varphi}
\oplus_{j=1}^{k+2} \OO_{\PP^2}(-b_j)\rightarrow
F_{HM}|_{\PP^2}\rightarrow 0.
\end{equation}
We prove that all of them are actually assumed by the HM-bundle on
some
plane $\PP^2 \subset \PP^4$.\\
\indent For each of these 6 resolutions, we then consider points of
the Grassmannian $\G(2,4)$ where this particular resolution is
assumed, forming the subvarieties that we call
$\VV_{(\underline{a},\underline{b})}$. Studying resolutions in
connection with jumping phenomena, we obtain a detailed description
of what
these subvarieties look like.\\
The results are summarized in the:\\
$\:$\\
\textbf{Main Theorem.} Let $\VV_{(\underline{a},\underline{b})}$ be
the subvarieties of the Grassmannian $\G(2,4)$ defined as follows:
$$
\VV_{(\underline{a},\underline{b})} := \{\PP^2 \subset \PP^4 \:|\: F_{HM}|_{\PP^2} \hbox{ has minimal free resolution
}(\underline{a},\underline{b})\}
$$
Then we have the filtration:
$$
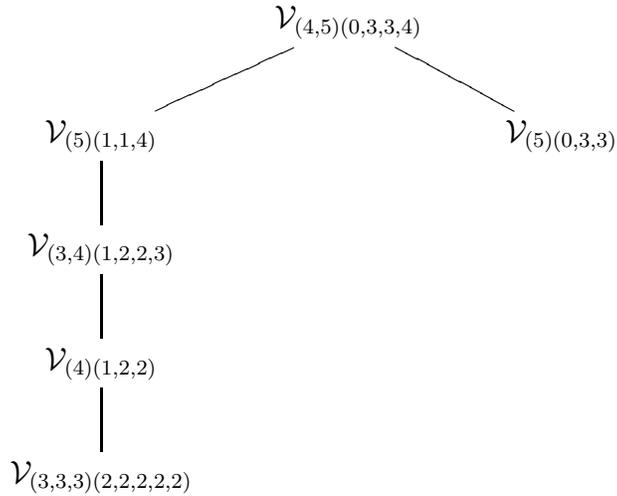
\xymatrix{&\VV_{(4,5)(0,3,3,4)}\ar@{-}[dl]\ar@{-}[dr]&\\
\VV_{(5)(1,1,4)}\ar@{-}[d]&&\VV_{(5)(0,3,3)}\\
\VV_{(3,4)(1,2,2,3)}\ar@{-}[d]&&\\
\VV_{(4)(1,2,2)}\ar@{-}[d]&&\\
\VV_{(3,3,3)(2,2,2,2,2)}&&}
$$
where: \begin{itemize}
       \item $\VV_{(5)(0,3,3)}$ is a smooth surface of degree 25, formed by jumping planes. It is the image of the well known Shioda's
       surface under a complete linear system;
       \item $\VV_{(4,5)(0,3,3,4)}$ consists of 25 smooth conics of jumping planes on the surface $\VV_{(5)(0,3,3)}$;
       \item $\VV_{(5)(1,1,4)}$ consists of 25 planes, formed by stable planes, each one containing one of the conics above described;
       \item $\VV_{(3,4)(1,2,2,3)}$ is an irreducible 4-fold, formed by stable planes;
       \item $\VV_{(4)(1,2,2)}$ is also formed by stable planes. It has dimension 5 and degree 5;
       \item $\VV_{(3,3,3)(2,2,2,2,2)}$ is an open subset of $\G(2,4)$, associated to the generic stable resolution.
     \end{itemize}
$\:$\\
\emph{Notation.} $V$ is a fixed 5-dimensional vector space over
$\C$, and $\PP^4=\PP(V)$  is the projective space of hyperplanes in
$V$, so that
$$
\HH^0(\OO_{\PP^4}(1))=V.
$$
$F_{HM}$ denotes the normalized Horrocks-Mumford bundle, with Chern
classes $c_1=-1$ and $c_2=4$. The stability of $F_{HM}$ (in the
Mumford-Takemoto sense) therefore simply means:
$$
\hh^0(F_{HM})= \dim (\HH^0(F_{HM}))=0
$$
as it is explained in \cite{OSS}, for example.\\

\indent I would like to thank Professor Hulek for the useful dicussion that we had in Ferrara and for the ideas
that led us to prove Corollary \ref{hulek}.\\
For the original idea of this paper, however, and for supervision and encouragement, I wish to express my deepest gratitude to Professor
Ottaviani.

\section{Minimal resolutions}
We are looking for all possible minimal free resolutions for
$F_{HM}|_{\PP^2}$. For this purpose, we use results on moduli spaces
contained in \cite{BS}, \cite{DM}, \cite{OV} and \cite{HL}.\\

\indent From Horrocks' Theorem (see \cite{Hor}) we learn that a
vector bundle on $\PP^n$ has homological dimension (\emph{i.e.} the
length of its minimal free resolution) at most $n-1$. This means
that we are actually looking for short exact sequences of the form:
\begin{equation}\label{risoluzione}
 0\rightarrow \oplus_{i=1}^{k} \OO_{\PP^2}(-a_i)\xrightarrow{\varphi}
\oplus_{j=1}^{k+2} \OO_{\PP^2}(-b_j)\rightarrow
F_{HM}|_{\PP^2}\rightarrow 0
\end{equation}
for some integers $a_i$, $b_j$ and $k$. Notice here that
$\rk(F_{HM})=2$.\\
 We assume that the two sequences of integers
$a_i$ and $b_j$ are indexed in nondecreasing order, and we call
$(\underline{a},\underline{b})=((a_1,...,a_k),(b_1,...,b_{k+2}))$
the \emph{associated pair}
to the resolution \ref{risoluzione}.\\
Our goal is now to determine bounds on $a_i$ and $b_j$, in order to
obtain all possible associated pairs for $F_{HM}|_{\PP^2}$. We have:
\begin{prop}\label{condizioni}
The integers $a_i$ and $b_j$ must satisfy the following six
conditions:
\begin{enumerate}
    \item${\sum_{i=1}^{k} a_i} - {\sum_{j=1}^{k+2} b_j} =
    -1;$
    \item${\sum_{i=1}^{k} {a_i}^2} - {\sum_{j=1}^{k+2} {b_j}^2} = 7;$
    \item$a_i\geq b_{i+2} + 1, \:\:\forall\: i=1...k;$
    \item$b_1\geq 0;$
    \item$2{b_{k+2}} + k \leq 11;$
    \item$2{{a_k} + k }\leq 12.$
    \end{enumerate}
\end{prop}
\begin{proof}
\emph{(1)} and \emph{(2)} come from asking the Chern classes of the bundle to be respectively  -1 and 4, as we must have:
\begin{align}
\nonumber \sum_{i=1}^k a_i -\sum_{j=1}^{k+2} b_j=c_1(F_{HM})=-1\\
\nonumber \sum_{i=1}^k a_i^2 -\sum_{j=1}^{k+2} b_j^2=2c_2(F_{HM})-c_1(F_{HM})^2=7.
\end{align}
Condition \emph{(3)} corresponds to the minimality of the resolutions (see \cite{BS}, Proposition 2).\\
Condition \emph{(4)} comes from a direct computation of the cohomology of the bundle restricted to a plane, which can be obtained from the
cohomology of the bundle on $\PP^4$, see \cite{HM}. More in detail, the Koszul complex splits in 2 short exact sequences:
\begin{equation}
\xymatrix@R-2ex@C-2ex{ 0\ar[r]&{\OO}_{{\PP}^4}(-2)\ar[r]&{{\OO}_{{\PP}^4}(-1)}^2\ar[dr]\ar[rr]
&&{\OO}_{{\PP}^4}\ar[r]&{\OO}_{{\PP}^2}\ar[r]&0\\
&&&\II_{{\PP^2},{\PP^4}}\ar[ur]\ar[dr]&&&\\
&&0\ar[ur]&&0&&\\
}
\end{equation}
which we can tensorize by $F_{HM}$ and thus obtain that for any $\PP^2 \subset \PP^4$, $\hh^0({F_{HM}|_{\PP^2}})=0$. This latter condition
implies
\emph{(4)}, as proved in \cite{DM}.\\
For \emph{(5)} and \emph{(6)} we modify an argument from \cite{OV} as follows. Using conditions \emph{(1)},
\emph{(2)} and \emph{(3)}, we have:
\begin{align}
\nonumber 7&= -{b_1}^2-{b_2}^2 + \sum_{i=1}^{k}({a_i}^2-{b_{i+2}}^2)=\\
\nonumber &=-{b_1}^2-{b_2}^2 +
\sum_{i=1}^{k}{(a_i-b_{i+2})(a_i+b_{i+2}-2b_2)}
+ 2b_2(b_1+b_2-1)\geq\\
\nonumber &\geq
-{b_1}^2-{b_2}^2+2\sum_{i=1}^{k}{(b_{i+2}-b_2)}+k+2b_2(b_1+b_2-1)\geq\\
\nonumber &\geq -{b_1}^2-{b_2}^2 + 2b_{k+2} - 2b_2 + k + 2b_2(b_1+b_2-1)=\\
\nonumber &=-[{b_1}^2+{b_2}^2-2b_2(b_1+b_2-2)] + 2b_{k+2} + k
\end{align}
Now looking at the hyperbola in the plane:
$$x^2-y^2-2xy+4y-4=0,$$
we get that, $\forall \:\:0\leq b_1\leq b_2 \in \Z$:
$${b_1}^2+{b_2}^2-2b_2(b_1+b_2-2)\leq 4.$$
Thus
$$
7 \geq -4 + 2b_{k+2} + k,
$$
which proves \emph{(5)}. Condition \emph{(6)} easily follows with an
identical argument as for \emph{(5)}.
\end{proof}

A brute force check after the conditions of Proposition
\ref{condizioni} proves the following:
\begin{cor}
The only pairs (thus possible resolutions) satisfying the conditions
of Proposition \ref{condizioni} are:
\begin{align}
\nonumber (5)\:&(0,3,3)\\
\nonumber (4,5)\:&(0,3,3,4)\\
\nonumber  \\
\nonumber (5)\:&(1,1,4)\\
\nonumber (3,4)\:&(1,2,2,3)\\
\nonumber (4)\:&(1,2,2)\\
\nonumber (3,3,3)\:&(2,2,2,2,2)
\end{align}
\end{cor}

\indent From a first look at the cohomology of bundles associated to
these six resolutions, we can immediately see that the first two are
not stable, whereas the other four do not have sections, \emph{i.e.} they are stable.\\
\indent It is also noteworthy that there are two couples of
resolutions that have the same cohomology: one is stable
($(4)(1,2,2)$ and $(3,4)(1,2,2,3)$), and one is not ($(5)(0,3,3)$
and $(4,5)(0,3,3,4)$). The reason why this can happen is that
minimality of the resolutions implies that every time we have a
constant map:
$$\OO(-a) \rightarrow \OO(-a),$$
it has to be the zero map (again, \cite{BS}, Prop. 2).\\
Anyhow, in case we wanted a cohomological distinction between those
resolution forming couples, we could obtain it by tensorizing for an
appropriate sheaf of forms.\\
\indent Notice that all bundles we are dealing with are
prioritaries, according to \cite{HL}.
\begin{rem}
In what follows we will often write $(\underline{a},\underline{b})$
meaning either the resolution, or the associated pair, or sometimes
even a bundle with that resolution (or better, its class of
isomorphism in the moduli space). What we are referring to will be
clear from the context.
\end{rem}

\section{Jumping phenomena}
It is known that $F_{HM}$ remains stable when it is restricted to
any hyperplane $\PP^3 \subset \PP^4$ (see \cite{BHM}, \cite{DS}, or
\cite{Hu}).\\
We give the following:
\begin{defn}
The subspace $X \subset \PP^4$ is a \emph{jumping space} for the HM-bundle iff $\HH^0(F_{HM}|_X) \neq 0$.
\end{defn}
Hence $F_{HM}$ has no jumping 3-spaces, but it admits jumping spaces
of lower dimension, planes and lines.\\
\indent For lines in particular we have a more subtle definition:
jumping lines are those where the bundle doesn't have the generic
splitting type entailed by Grauert-M\"{u}lich Theorem. The measure
of how much a jumping line drifts away from the generic case is given by its jumping order.\\
In other words:
\begin{defn}
A line $\ell \subset \PP^4$ is a \emph{jumping line of order $k$} for HM (also called a \emph{$k$-jumping line}) iff $\HH^0(F_{HM}|_{\ell}(-k))
\neq 0$.
\end{defn}
HM-bundle's jumping lines and planes have been deeply
studied by Barth, Hulek and Moore in \cite{BHM}.\\
The structure of jumping planes is explained in the following
result:
\begin{thm}\label{pianisalto}
A plane $\PP^2 \in \G(2,4) \subset \PP^9$, with Pl\"{u}ckerian
coordinates $p_{ij}$, is a jumping plane iff $\rk(A \otimes
\OO_{\PP^2}) \leq 1$, where
\begin{equation}\label{A}
A \otimes \OO_{\PP^2}=\left(%
\begin{array}{ccccc}
  p_{23} & p_{34} & -p_{04} & p_{01} & p_{12} \\
  p_{14} & -p_{02} & -p_{13} & -p_{24} & p_{03} \\
\end{array}%
\right).
\end{equation}
Jumping planes thus form a smooth surface of degree 25 in $\PP^9$,
called $S_{25}$, which is birational to Shioda's modular surface
$S(5)$. $S_{25}$ may also be seen as the transverse intersection of
the Grassmannian $\G(2,4)$ with the (suitably normalized) Segre
variety $\PP^1 \times \PP^4$.
\end{thm}
\begin{proof}
The result simply derives from translating in Pl\"{u}cker
coordinates the condition $\hh^0 \neq 0$.  Note that the matrix $A$
is the one used to define the monad in (\ref{monadeHM}).\\
Let our $\PP^2 \subset \PP^4$ be defined as zero locus of the 2
hyperplanes $\xi^*$ and $\zeta^*$. Then
$\PP^2=\{\sum_0^4\xi_ix_i=\sum_0^4\zeta_ix_i=0\}$. Restricting $A$
to $\PP^2$ means now contracting the 2-forms $a_{ij}$ with the
2-form $\xi^* \wedge \zeta^*$ associated to the plane, i.e. $ A
\otimes \OO_{\PP^2}=(\xi^* \wedge \zeta^*(a_{ij}))$, hence
Pl\"{u}ckerian
coordinates $p_{ij}=\xi_i\zeta_j-\xi_j\zeta_i$ appear.\\
For details, see \cite{DS}, Proposition 2, and \cite{BHM}, section
3.4.
\end{proof}

\indent Notice that we have a natural map $S_{25} \rightarrow
\PP^1$, and that the general fibers are transverse linear sections
of $\G(2,4)$, hence they are elliptic normal quintics in $\PP^4$.\\

\indent For the central role of Shioda's modular surface in the
understanding of HM-bundle jumping phenomena, we refer the reader to
the detailed
articles \cite{BHM} and \cite{BH}, and to \cite{Hu}.\\
Here we just limit ourselves to mention the result on jumping lines,
which we will need later on. Our bundle admits jumping lines of
order 1, 2 and 3. More precisely:
\begin{thm}[{[BHM], section 1.2 and 5.5}]\label{j1j2j3}
Let $J_i(F_{HM}) \subset \G(1,4)$, for $i=1,2,3$, be the
subvarieties of jumping lines for the HM-bundle.\\
$J_1(F_{HM})$, is a rational, irreducible 4-fold. It is smooth
outside of $J_2(F_{HM})$. The variety $J_2 \subset J_1$ has
dimension 2, and it is birational to Shioda surface $S(5)$. In fact
it is nothing but $S(5)$ with its 25 sections blown down to 25
singular points, forming $J_3(F_{HM})$. $J_3$ is formed by 25
points, corresponding in $\PP^4$ to 25 skew lines, the so called
Horrocks-Mumford lines.
\end{thm}

\indent From this first look at jumping spaces, we can already infer
that at least two out of the six possible resolutions are actually
assumed
by $F_{HM}|_{\PP^2}$: one of them will be stable, the other not.\\
The situation so far is then:
\begin{equation}\label{situazione}
S_{25} \subset \G(2,4) \subset \PP^9
\end{equation}
where $S_{25}$ ``contains'' the 2 unstable resolutions, the other
4 being ``contained'' in $\G(2,4) \setminus S_{25}$.\vspace{0.5cm}\\
\indent Keeping in mind Horrocks-Mumford bundle's jumping phenomena,
we now
want to analyze jumping lines admitted by each of our six resolutions.\\
The situation is described in the following result:
\begin{prop}\label{risoluzionirettesalto}$\:$
\begin{itemize}
\item Part A
\begin{itemize}
\item$(5)(0,3,3)$ admits jumping lines of order 1 and 2.\\
Take the 4 points of intersection of the 2 conics $q_i$, zero loci
of $q_i: \OO_{\PP^2}(-5) \rightarrow \OO_{\PP^2}(-3)$. The generic
structure then is the following: the 4 pencils of lines through
these 4 points are the jumping lines, while the 6 lines
connecting pairs among the 4 points are double jumping lines.\\
The 4 points are exactly those where the unique section violating
stability vanishes.
\item$(4,5)(0,3,3,4)$ generically admits only 1 jumping line $\ell_1$, of order 3.
It is is the zero locus of $\ell_1: \OO_{\PP^2}(-5) \rightarrow
\OO_{\PP^2}(-4)$.
\item$(5)(1,1,4)$ generically admits only 1 triple jumping line $\ell_2$. It
is the zero locus of $\ell_2: \OO_{\PP^2}(-5) \rightarrow
\OO_{\PP^2}(-4)$.
\item$(3,4)(1,2,2,3)$ generically admits only 1 double jumping line $\ell_3$, zero locus
 of $\ell_3: \OO_{\PP^2}(-4) \rightarrow \OO_{\PP^2}(-3)$.
\item$(4)(1,2,2)$ admits only jumping lines of order 1.\\
Those are, in the generic case, the 6 lines connecting pairs among
the 4 points of intersection of the conics zero loci of $q_i:
\OO_{\PP^2}(-4) \rightarrow \OO_{\PP^2}(-2)$
\end{itemize}
\item Part B
\begin{itemize}
\item$(3,3,3)(2,2,2,2,2)$ admits only 1-jumping lines.\\
They can be either 6 lines with normal crossing, or they can form a
conic.
\end{itemize}
\end{itemize}
\end{prop}

\begin{proof}
\emph{of Part A.}\\
We will not go through all the 6 cases in detail, because we use
exactly the same argument for the 5 resolutions of Part A.\\
\indent Take for example a vector bundle $E$ associated with
$(5)(0,3,3)$. We have
\begin{equation}\label{(5)(033)}
0 \rightarrow \OO_{{\PP}^2}(-5) \xrightarrow{\varphi}
2{\OO_{{\PP}^2}(-3)} \oplus \OO_{{\PP}^2} \rightarrow E \rightarrow
0
\end{equation}
where
$$
^t\varphi=\begin{tabular}{|c|c|c|}
  \hline
  $q_1$ & $q_2$ & $r$ \\
  \hline
\end{tabular}
$$
and the $q_i$'s are quadratic polynomials, while $r$ has degree
5.\\Once dualized and restricted to a generic line $\ell \subset
\PP^4$, the sequence (\ref{(5)(033)}) becomes (remember that
$E^*=E(-c_1)=E(1)$):
\begin{equation}\label{succ da ritwistare}
0 \rightarrow E|_{\ell} \rightarrow \OO_{{\PP}^1}(-1) \oplus
2{\OO_{{\PP}^1}(2)} \rightarrow \OO_{{\PP}^1}(4) \rightarrow 0
\end{equation}
If we look at its cohomology sequence, we get:
\begin{align}
\nonumber 0 \rightarrow \HH^0(E|_{\ell}(-1)) \rightarrow 2V
&\xrightarrow{\Phi} S^3(V)\\
\nonumber  (\LL_1, \LL_2) &\mapsto (\LL_1 \cdot {q_1}|_{\ell} + \LL_2 \cdot {q_2}|_{\ell}).
\end{align}
Now the two schemes ${q_1}|_{\ell}$ and ${q_2}|_{\ell}$ both consist
of two points, and there is a nonzero $(\LL_1, \LL_2) \in \ker\Phi$
iff there is a common point between these two schemes, that is iff
$\ell\in q_1\cap q_2$\\
In other words, when $\LL_1$ and $\LL_2$ vary, $\ell$ describes four
pencils of lines
through the 4 points of intersections of the two conics.\\
Moreover these 4 points coincide with the zero locus of the generic
section $s$ that breaks the stability.\\
To see this, first remark that $\hh^0(E)=1$ implies that we have
only one section, and that $\hh^0(E(-1))=0$ and $c_2=4$ imply
respectively that the zero locus has dimension 0 and degree 4. Now
let
\begin{equation}\label{sezione(5)(033)}
0 \rightarrow \OO_{{\PP}^2} \xrightarrow{s} E \rightarrow \II_Z(-1)
\rightarrow 0
\end{equation}
be the short exact sequence of defining the subscheme
$Z=\{P_1,...,P_4\}$, of points where $s$ vanishes.\\
Putting (\ref{sezione(5)(033)}) together with (\ref{(5)(033)}), we
obtain the diagram:
$$
\xymatrix{&&&0\ar[d]&\\
          &&&\OO_{{\PP}^2}\ar[d]^s&\\
          0\ar[r]&\OO_{{\PP}^2}(-5)\ar[r]^{(q_1,q_2,r)}&2{\OO_{{\PP}^2}(-3)} \oplus
          \OO_{{\PP}^2}\ar[r]&E\ar[d]\ar[r]&0\\
          &&&\II_Z(-1)\ar[d]&\\
          &&&0&
 }
$$
Now, completing the diagram, we get our thesis:
$$
\xymatrix{&&0\ar[d]&0\ar[d]&\\
          &&\OO_{{\PP}^2}\ar[d]\ar@{=}[r]&\OO_{{\PP}^2}\ar[d]^s\ar@{-->}[dl]&\\
          0\ar[r]&\OO_{{\PP}^2}(-5)\ar@{=}[d]\ar[r]^{(q_1,q_2,r)}&2{\OO_{{\PP}^2}(-3)}
          \oplus \OO_{{\PP}^2}\ar[r]\ar[d]&E\ar[d]\ar[r]&0\\
          0\ar[r]&\OO_{{\PP}^2}(-5)\ar[r]^{(q_1,q_2)}&2{\OO_{{\PP}^2}(-3)}\ar[d]\ar[r]
          &\II_Z(-1)\ar[d]\ar[r]&0\\
          &&0&0&
 }
$$
To prove the part about the 6 double jumping lines, we repeat the previous argument, except that this time we twist (\ref{succ da ritwistare})
by $\OO_{\PP^1}(-2)$ and the cohomology sequence looks:
\begin{align}
\nonumber 0 \rightarrow \HH^0(E|_{\ell}(-2)) \rightarrow 2{\C} &\xrightarrow{\Phi} S^2(V)\\
\nonumber  (\alpha, \beta) &\mapsto (\alpha \cdot {q_1}|_{\ell} +
\beta \cdot {q_2}|_{\ell})
\end{align}
It is clear that, generically, to make the map $\Phi$ vanish we
need:
$$\alpha \cdot q_1 |_{\ell}= -\beta \cdot q_2|_{\ell}.$$
\noindent In other words we want the 4 points given on $\ell$ by the
2 conics to coincide 2 by 2. Generically this is possible for
appropriate values of $\alpha$ and $\beta$, iff $\ell$ is exactly
one of the 6
lines passing through the 4 intersections of $q_1$ and $q_2$, that is what we wanted.\\
\\
\emph{Proof of Part B.}\\
We explain the generic resolution separately, because we treated
this last case with different tools, since we deal with a Steiner
bundle.
We use results from \cite{AO} and \cite{Va}.\\
\indent From the cohomology table, we know that only order 1 jumping
lines
are admitted. They can be either 6 lines with normal crossing, or they can form a conic in our plane.\\
First, let $W(S)$ be the scheme (see \cite{AO}) of unstable
hyperlanes of a Steiner bundle $S \in S_{n,k}$ on $\PP^n$, with dual
resolution
$$
0 \rightarrow S^* \rightarrow W \otimes \OO_{{\PP}^n}
\xrightarrow{f_A} I \otimes \OO_{{\PP}^n}(1) \rightarrow 0,
$$
where $W$ and $I$ are complex vector spaces of dimension $n+k$ and
$k$ respectively. Then we have:
\begin{prop}[{\cite{AO}, section 3.8}] Let $p_V$ be the projection of the
Segre variety $\PP(V) \times \PP(I)$ on the $\PP(V)$. Then
$$
{W(S)}_{red}=p_V{[\PP(W) \cap (\PP(V) \otimes \PP(I))]_{red}}
$$
\end{prop}
In our case where $(n,k)=(2,3)$, the generic Steiner bundle is
logarithmic. In fact, the generic $\PP^4$ linearly embedded in
$\PP^8$ meets the Segre variety $\PP^2 \otimes \PP^2$ in $\deg(\PP^2
\otimes
\PP^2)=6=n+k+1$ points.\\
Contrary to what happens for the resolution $(4)(1,2,2)$, the six
lines don't have a particular configuration, but normal crossing, as
it is stated in \cite{AO}, Theorem 5.6.\\
In case the jumping lines are not 6, they can only be infinite and thus form a conic (see \cite{AO}). In particular, we get that if a generic
linear $\PP(W)$ meets $\PP(V) \otimes \PP(I)$ in $n+k+2$ points, then it meets it in infinitely many points.
\end{proof}
\indent Now we would like to improve (\ref{situazione}), and to obtain a more detailed filtration of the Grassmannian.\\
If we look at Theorem \ref{pianisalto}, we see that the
characterization of jumping planes has been obtained by translating
in Pl\"{u}cker coordinates a
cohomological condition.\\
Can we repeat this argument?\\
Looking at the cohomology table of the bundles associated to the 6
resolutions:
$$
\begin{tabular}{r||c|c|c|c}
  & $\hh^0(E(3))$ & $\hh^0(E(2))$ & $\hh^0(E(1))$ & $\hh^0(E)$ \\
  \hline \hline
   (4,5)(0,3,3,4)& 12 &6 &3  &1  \\
  \hline
  (5)(0,3,3)& 12 & 6 & 3 & 1 \\
  \hline
   (5)(1,1,4)& 12 & 6  & 2 & 0 \\
  \hline
   (3,4)(1,2,2,3)& 12 & 5 & 1 & 0 \\
  \hline
   (4)(1,2,2)& 12 & 5 & 1 & 0 \\
  \hline
   (3,3,3)(2,2,2,2,2)& 12 & 5 & 0 & 0 \\
\end{tabular}
$$
it is clear that the best condition to translate would be $\hh^0(E(1))$.\\
Doing this is not so easy as it has been for the previous twist,
though. Take the display associated to the monad (\ref{monadeHM}):
$$
\xymatrix{&&0\ar[d]&0\ar[d]&\\
0\ar[r]&5{\OO_{\PP^4}(-1)}\ar[r]\ar@{=}[d]&\Ker(A)\ar[r]\ar[d]&F_{HM}\ar[d]\ar[r]&0\\
0\ar[r]&5{\OO_{\PP^4}(-1)}\ar[r]^B&2{\Omega}^2_{{\PP}^4}(2)\ar[r]\ar[d]^A&
\coker(B)\ar[r]\ar[d]&0\\
&&5{\OO_{\PP^4}}\ar@{=}[r]\ar[d]&5{\OO_{\PP^4}}\ar[d]&\\
&&0&0&
 }
$$
Looking at it, from the first line we get:
$$
\hh^0({F_{HM}(1)}|_{{\PP}^2}) = \hh^0({\Ker(A(1))}|_{{\PP}^2})-5
$$
and from the first row:
$$
\hh^0({\Ker(A(1))}|_{{\PP}^2})= 20 - \rk(A(1) \otimes \OO_{{\PP}^2})
$$
Unlike the previous twist (see \cite{DS}, Prop 2 for details), this
time we are not able to determine explicitly the map $A(1) \otimes
\OO_{{\PP}^2}$. Nevertheless, we can compute
its rank by embedding everything in a bigger space.\\
Let's see this in detail. We have an isomorphism:
\begin{equation}\label{isomorfismoomega}
\HH^0({\Omega}^2_{{\PP}^4}(3)) \cong \HH^0({{\Omega}^2_{{\PP}^4}(3)}|_{{\PP}^2}).
\end{equation}
This can be verified by tensorizing short exact sequence of definition of an hyperplane $\PP^3 \subset \PP^4$ by the sheaf
${\Omega}^2_{\PP^4}(3)$:
$$
0 \rightarrow {\Omega}^2_{{\PP}^4}(2) \rightarrow {\Omega}^2_{{\PP}^4}(3) \rightarrow {{\Omega}^2_{{\PP}^4}(3)}|_{{\PP}^3} \rightarrow 0.
$$
Looking at the cohomology, since ${\Omega}^2_{{\PP}^4}(2)$ has vanishing cohomology, we get that:
\begin{equation}
\HH^0({\Omega}^2_{{\PP}^4}(3)) \cong \HH^0({{\Omega}^2_{{\PP}^4}(3)}|_{{\PP}^3}).
\end{equation}
Now since ${\Omega}^2_{{\PP}^4}(2)|_{\PP^3}$ has vanishing cohomology as well, repeating the argument for an (hyper)plane $\PP^2 \subset \PP^3$
we get the desired
isomorphism (\ref{isomorfismoomega}).\\
Now we use the Pl\"{u}cker embedding:
\begin{align}
\nonumber V &\xrightarrow{\psi(\omega)} \bigwedge^3V\\
\nonumber v &\mapsto v \wedge \omega
\end{align}
where $\psi$ is given by the 2-form $\omega=x^* \wedge y^*$, once we
have defined our plane $\PP^2 \subset \PP^4$ as the zero locus of
the two hyperplanes $x^*$ and $y^*$:
$$
\PP^2=\{x^*=y^*=0\}
$$
In this case $W=\im(\psi(\omega))$, where $\PP^2=\PP(W)$.\\
Thus we have:
\begin{equation}\label{diagrammafinale}
\xymatrix{2{\HH^0({\Omega}^2_{\PP^4}(3)|_{\PP^2})}\ar[r]^{A(1)
\otimes \OO_{\PP^2}} & 5{\HH^0(\OO_{\PP^2}(1))}\ar@{^{(}->}[r]^{5i}&
5 \bigwedge^3
\HH^0(\OO_{\PP^4}(1))\\
2\HH^0({\Omega}^2_{\PP^4}(3))\ar[u]^{\thicksim}\ar[ur]\ar[r]^{A(1)}
&5{\HH^0(\OO_{\PP^4}(1))} \ar[u]_{5\pi}\ar[ur]_{5\psi(\omega)}& }
\end{equation}
All in all, we have obtained:
\begin{equation}\label{defM}
\rk(A(1) \otimes \OO_{{\PP}^2})=\rk(5\psi(\omega) \circ A(1))
\end{equation}
We are now ready to state the first result.
\begin{prop}\label{caratterizzcoomologica}
Let $\pi$ be a plane in $\G(2,4) \subset \PP^9$. Let
$$
M=A(1) \otimes \OO_{{\PP}^2}: 2\HH^0({\Omega}^2_{\PP^4}(3))
\rightarrow 5 \wedge^3 \HH^0(\OO_{\PP^4}(1))
$$
be the $20 \times 50$ matrix constructed above. The rank of $M$
determines the minimal resolution of $F_{HM}|_{\pi}$. More
precisely, we have:
\begin{itemize}
    \item $\rk M=15$ $\Leftrightarrow$ $F_{HM}|_{\pi}$ has
    resolution $(3,3,3)(2,2,2,2,2)$;
    \item $\rk M=14$ $\Leftrightarrow$
    $F_{HM}|_{\pi}$ has resolution either $(4)(1,2,2)$ or \\$(3,4)(1,2,2,3)$;
    \item $\rk M=13$ $\Leftrightarrow$
    $F_{HM}|_{\pi}$ has resolution $(5)(1,1,4)$;
    \item $\rk M=12$ $\Leftrightarrow$
    $F_{HM}|_{\pi}$ has resolution either $(5)(0,3,3)$ or \\$(4,5)(0,3,3,4)$.
\end{itemize}
\end{prop}
\begin{rem}Clearly this classification, being based on a cohomological
criterion, doesn't take into account the differences between those
couples of resolutions with the same cohomology (see paragraph 2).
\end{rem}

\section{The $M_k$ subvarieties}
Let's define the subvarieties of the Grassmannian $\G(2,4)$:
$$
M_{k} := \{\pi \in \G(2,4)\:|\:\rk M(\pi) \leq k \}.
$$
\noindent We want to know more about $M_{12}$,
$M_{13}$, $M_{14}$ and $M_{15}$.\\
Notice that even though we do have explicit equations for them,
these equations are almost impossible to handle, even with the aid
of a computer.\\
\begin{rem} All the four $M_k$'s are non empty.
In fact we are able to find out points where $\rk(M)$ assumes all its four possible values (see next Proposition \ref{dimlocali} for
details).\\
This already implies that we can confirm the presence of the generic
resolution $(3,3,3)(2,2,2,2,2)$ between those really assumed by $F_{HM}|_{\PP^2}$.\\
\end{rem}
We can make some other useful observations. First $M_{15}$, being
associated to the generic resolution, hence to the open condition of
maximum rank, is an open subset of $\G(2,4)$, so $\dim M_{15}=\dim
\G(2,4)=6$. Moreover we know from the previous remark that $M_{15}
\neq
 \emptyset$. From now on, recalling the notation given in the Introduction, we
will call this last one $\VV_{(3,3,3)(2,2,2,2,2)}$.\\
Second, the
subvariety $M_{12}$ is nothing else than $S_{25}$, the surface
in $\PP^9$ that parametrizes jumping planes.\\
So now the situation described in (\ref{situazione}) has been
improved quite a lot:
\begin{equation}\label{nuovasituazione}
S_{25}=M_{12} \subset M_{13} \subset M_{14} \subset M_{15}=\VV_{(3,3,3)(2,2,2,2,2)} \subset \G(2,4) \subset \PP^9.
\end{equation}
\vspace{0.3cm}\\ \indent Another useful step in the direction of
understanding the $M_k$'s
structure is computing their dimension.\\
\indent Let's take for example $M_{13}$. The matrix $M$ clearly
defines a map $\vartheta$:
$$
\xymatrix@C-2ex{\G(2,4) \ar[rr] \ar@{^{(}->}[dr]&&\PP(\mathbb{M}(20\times50))\\
&\PP^9 \ar[ur]_{\vartheta}&\\ }
$$
where $\vartheta$ is linear in the coordinates $p_{ij}$.\\
Now let $R_{13} \subset \PP(\mathbb{M}(20\times50))$ the subspace
consisting of matrices whose rank is 13. We have:
$$
M_{13}= {\vartheta}^{-1}(R_{13}) \cap \G(2,4).
$$
For each point (plane) $\pi \in \G(2,4)$ we could then obtain the
tangent space $\T_{\pi}(M_{13})$, simply by deriving the equations
given by the appropriate minors of $M$, together with the 5
Pl\"{u}cker quadrics, of course.\\
Obviously these equations, the $14 \times 14$ minors of a $20 \times
50$ matrix, are too heavy to handle. But if we are satisfied with
local dimension (\emph{i.e.} local tangent space), we can use the
point where we know that $\rk(M)=13$.\\
In other words, we have a
$13 \times 13$ minor $\mu$ which we know not to be zero. Then we can
compute only those minors $14 \times 14$ obtained by adding 1 row
and 1 column to $\mu$. This diminishes considerably the
computational cost
of the operation, and can be done without great efforts.\\
\indent Thanks to these simple observations, and with the aid of Macaulay2 computer system \cite{Mac2}, we obtain:
\begin{prop}\label{dimlocali}
The subvarieties $M_k$ have dimension:
\begin{itemize}
    \item $\dim M_{12}=2$, locally around $\pi=(1,0,0,0,0,0,0,0,0,0)$;
    \item $\dim M_{13}=2$, l.a. $\rho=(1,1,1,1,0,0,0,0,0,0)$;
    \item $\dim M_{14}=5$, l.a. $\sigma=(1,1,0,0,1,0,0,0,0,0)$.
\end{itemize}
\end{prop}
Remark that $\dim M_{12}$ is exactly the expected one, since
$M_{12}$ is the surface $S_{25}$.\\

\indent Going further on with the study of the $M_k$ subvarieties,
we use Proposition \ref{risoluzionirettesalto}. Let's take a jumping
line $L_{k,j}$ of order 3 for which we do have explicit equations,
from \cite{BHM}:
$$
L_{k,j}=\{z_k=z_{k+2}+e^{j\frac{2\pi\imath}{5}}z_{k+3}=z_{k+1}+e^{3j\frac{2\pi\imath}{5}}z_{k+4}=0\},
\:\: k,j=0...4.
$$
The incidence variety:
$$
\{\pi \in \G(2,4) \:|\: \pi \supset L_{k,j}\}
$$
forms a plane $\PP^2 \subset \PP^9$. This plane cuts the surface
$S_{25}$ in a smooth conic.\\
If we recall that $L_{k,j}$ may be contained only in a plane with
resolution either $(5)(1,1,4)$ if it is stable, or $(4,5)(0,3,3,4)$
if it is not, we easily have:

\begin{thm}\label{3-salto}
Each of the 25 (skew) jumping lines of order 3 is contained in a
smooth conic $C \subset \PP^2$ of jumping planes, on which the
restriction of $F_{HM}$ takes resolution $(4,5)(0,3,3,4)$.\\
For each of these lines then there exists a set $\{\PP^2 \setminus
C\}$ of stable planes containing it, with resolution $(5)(1,1,4)$.
\end{thm}
\begin{proof}
Now we explicitly use the group of symmetries of the HM bundle, for
which we refer the reader to \cite{HM}, \cite{Hu} and \cite{Ma}. In
fact the property we are checking is invariant under the action of
the symmetry group of $F_{HM}$, and this group acts transitively on
the 25 Horrocks-Mumford lines (see \cite{HM}). This means that what
we prove for one $L_{k,j}$ is valid for all
the 25 triple jumping lines.\\
Take for example the line:
$$L_{0,0}=\{z_0=z_2+z_3=z_1+z_4=0\}.$$
A plane $\pi=\{{\sum_{i=0}^{4}f_iz_i}={\sum_{i=0}^{4}g_iz_i}=0\}$
will contain $L_{0,0}$ iff:
\begin{equation}\label{richiesta}
\rk \left(%
\begin{array}{ccccc}
  f_0 & f_1 & f_2 & f_3 & f_4 \\
  g_0 & g_1 & g_2 & g_3 & g_4 \\
  1 & 0 & 0 & 0 & 0 \\
  0 & 0 & 1 & 1 & 0 \\
  0 & 1 & 0 & 0 & 1 \\
\end{array}%
\right) \leq 3.
\end{equation}
Translating condition (\ref{richiesta}) in Pl\"{u}cker coordinates
$p_{ij}=f_ig_j-f_jg_i$, we obtain:
$$
 \left\{
\begin{array}{l}
p_{14}=p_{23}=0\\
p_{01}=p_{04}=a\\
p_{02}=p_{03}=b\\
p_{13}=p_{12}=-p_{24}=-p_{34}=c\\
\end{array}
\right.
$$
Now if we ask for the plane $\pi$ to be unstable, hence to satisfy
Theorem \ref{pianisalto}, we finally get the smooth conic $ab=c^2$.
\end{proof}

\begin{cor}\label{hulek}
The surface $S_{25}$ contains exactly 25 conics, which are the image
of the 25 sections of the Shioda's surface.
\end{cor}
\begin{proof}
The 25 smooth conics of Theorem \ref{3-salto} are exactly the image
of the 25 sections of the Shioda's modular surface $S(5)$ under the
linear system
$$
S(5) \xrightarrow{|I+3F|} S_{25},
$$
where $I$ and $F$ are classes of divisors of $S(5)$ (see \cite{BH} and
\cite{BHM} for details).\\
This is because if we take any curve $C$ that doesn't correspond
neither to a section nor to a singular fibre, then $C$ must have
degree $\geq 3$, since both $C.I$ and $C.F$ are strictly greater
than zero (again, we refer to \cite{BH} and \cite{BHM}). Thus the
only curves with degree 2 are the 25 sections.
\end{proof}

Theorem \ref{3-salto} confirms the presence of $(4,5)(0,3,3,4)$ and
 $(5)(1,1,4)$, because we need to have both stable and unstable planes
 to contain jumping lines of order 3.\\
 \indent If we analyze this result from the $M_k$'s point of
view, we can distinguish inside $M_{12}=S_{25}$ two subvarieties:
\begin{align}
\nonumber &\VV_{(4,5)(0,3,3,4)} \leftrightarrow (4,5)(0,3,3,4)\\
\nonumber &\VV_{(5)(0,3,3)} \leftrightarrow (5)(0,3,3),
\end{align}
where of course $ M_{12}= \VV_{(4,5)(0,3,3,4)} \cup \VV_{(5)(0,3,3)}= S_{25} $.\\
Then Theorem \ref{3-salto} takes us to this
\begin{cor}$\:$
\begin{itemize}
    \item $\VV_{(4,5)(0,3,3,4)}=$ 25 smooth conics;
    \item $\VV_{(5)(0,3,3)}= S_{25} \setminus \VV_{(4,5)(0,3,3,4)}$;
    \item $M_{13} =:\VV_{(5)(1,1,4)} =$ 25 planes.
\end{itemize}
\end{cor}
Remark that $\VV_{(5)(1,1,4)}$ is reducible.\\
This result makes us hope we could obtain information on $M_{14}$
repeating a similar argument for a double jumping line.\\
Just like we have done for $M_{12}$, we make the distinction:
\begin{align}
\nonumber &\VV_{(3,4)(1,2,2,3)} \leftrightarrow (3,4)(1,2,2,3)\\
\nonumber &\VV_{(4)(1,2,2)} \leftrightarrow (4)(1,2,2).
\end{align}
Now recall that the resolution $(3,4)(1,2,2,3)$ is the only stable
one that admits jumping lines of order 2. We have that:
\begin{prop}$\VV_{(3,4)(1,2,2,3)}$ is irreducible of dimension 4.
\end{prop}
\begin{proof}
 Take a 2-jumping line, and recall that 2-jumping lines form a surface $J_2
\subset \G(1,4)$ which is smooth outside 25 points. Since (once
again from Proposition \ref{risoluzionirettesalto}) $(3,4)(1,2,2,3)$
admits only 1 jumping line of order 2, we have the fibration:
$$
\xymatrix{\VV_{(3,4)(1,2,2,3)} \ar@{>>}[d]^g& \pi \ar@{|-_{>}}[d]\\
J_2 & \ell}
$$
that associates to each plane $\pi \in \VV_{(3,4)(1,2,2,3)}$ the
only 2-jumping line $\ell$ contained in it. Each fiber has dimension
2. Moreover, $J_2$ is irreducible, and for every $\ell \in J_2$,
$g^{-1}(\ell)\simeq \PP^2$ is irreducible. This concludes the proof.
\end{proof}

\indent From \cite{BHM} we learn that generically each double
jumping line is contained
in exactly 3 jumping planes.\\
As a first consequence, now we can also confirm the presence of
resolutions $(3,4)(1,2,2,3)$ and $(5)(0,3,3)$,
for we must have stable and unstable planes containing jumping lines of order 2.\\
Secondly, with the notation used above, for $\ell \in J_2$ generic:
\begin{equation}\label{3punti}
g^{-1}(\ell)\cap S_{25}\:= 3 \points.
\end{equation}
As already remarked, the above described property is generic. In
fact when working out some examples explicitly with the aid of a
computer, sometimes a degenerate case comes out, where instead of 3
points, the intersection (\ref{3punti}) consists of a line plus an
external point.\\
\indent The observation of this phenomenon takes us to stress out
the fact that
there are lines entirely contained in $S_{25}$.\\
This we already knew from \cite{BHM}, section 3.4. In $ S_{25}$ we
find 12 pentagons, that are the 12 singular fibres of the Shioda's
surface, isomorphically mapped in $S_{25}$. It means we have at
least $12 \times 5=60$ lines inside our surface.\\
We can show that:
\begin{cor}\label{60rette}
The surface $S_{25}$ contains exactly 60 lines.
\end{cor}
\begin{proof}
We use the same argument that entailed Corollary \ref{hulek}. If
there is a curve in $S_{25}$ with degree 1, then it must be
contained in a singular fibre, otherwise it would have degree $\geq
3$.
\end{proof}

\indent As our last step, we want to study the subvariety
$\VV_{(4)(1,2,2)}$.\\
Of course we deal with the usual computing problems. Unfortunately
in this case jumping phenomena do not help, since we have another
stable resolution that admits jumping lines of order 1,
the generic $(3,3,3)(2,2,2,2,2)$.\\
Yet we have explicitly found out a point
$$\widetilde{\pi}=(1,1,1,0,0,0,0,0,0,0) \in \G(2,4)$$
where the resolution $(4)(1,2,2)$ is assumed, so that our subvariety $\VV_{(4)(1,2,2)}$
 is non-empty.\\
From Proposition \ref{dimlocali}, we know that (locally) it has
dimension 5. We want to compute its degree. The result we obtain is:
\begin{prop}\label{gradoM14g}
The subvariety formed by stable planes $\VV_{(4)(1,2,2)}$ is a
5-fold of degree 5.
\end{prop}
\begin{proof}
$\VV_{(4)(1,2,2)} \subset \G(2,4)$ has codimension 1, a hypersurface
inside the Grassmannian. Now let's take a line $\PP^1 \subset
\G(2,4)$, paramatrized by the coordinates $(s,t)$. The intersection
$\VV_{(4)(1,2,2)} \cap \PP^1$ is formed by a number of points equal
to the degree we are looking for:
$$
\deg (\VV_{(4)(1,2,2)} \cap \PP^1) = \deg (\VV_{(4)(1,2,2)}),
$$
but this time we are dealing with a Principal Ideal Domain, which
means that $I(\VV_{(4)(1,2,2)} \cap \PP^1)=(f)$, where $f \in
\C[s,t]$ is an homogeneous polynomial of degree = $\deg
(\VV_{(4)(1,2,2)})$.\\
The Greatest Common Divisor of two (suitably choosen) $15 \times 15$
minors is a degree 5 polynomial in $\C[s,t]$, thus:
$$
\deg (\VV_{(4)(1,2,2)}) \leq 5.
$$
In fact, the degree turns out to be exactly 5.\\
To show this, we use Invariant Theory. Let's take the standard
action of the Heisenberg group $H_5$ on $V$ (see \cite{HM} and
\cite{Ma}), and the induced action on $\wedge^2 V$. From the
computations contained in the appendix of Manolache's work, it is
easy to see that the symmetric powers
$$
\s^i(\wedge^2 V)
$$
do not contain trivial summands for $1 \leq i \leq 4$, so that the
degree of our subvariety must be 5.
\end{proof}

\begin{rem}
From the proof of Proposition \ref{gradoM14g} we can also infer that
the hypersurface we have found inside $\G(2,4)$ is irreducible.\\
Unfortunately we are not yet able to exclude the presence of some
components of lower dimension in $\VV_{(4)(1,2,2)}$.
\end{rem}


\begin{thebibliography}{50}

\bibitem{AO}V.Ancona and G.Ottaviani, \emph{Unstable hyperplanes for Steiner
bundles and multidimensional matrices}, Adv. Geom. 1(2), 2001,
165-192.

\bibitem{Ba}W.Barth, \emph{Some properties of stable rank 2 vector
bundles on $\PP^n$}, Math.Ann. 226 (1977), 125-150

\bibitem{BH}W.Barth and K.Hulek, \emph{Projective modules of Shioda modular
surfaces}, Manuscripta math. 50 (1985), 73-132.

\bibitem{BHM}W.Barth, K.Hulek and  R.Moore, \emph{Shioda's modular surface $S(5)$
and the Horrocks-Mumford bundle}, in \emph{Vector bundles on
algebraic varieties. Papers presented at the Bombay colloquium
1984}, 35-106, Oxford University Press, Bombay 1987.

\bibitem{BS}G.Bohnhorst and H.Spindler, \emph{The stability of certain vector bundles
on $\PP^n$}, Lecture Notes in Mathematics (1992), Springer, 1507,
39-50.

\bibitem{Co}A.Comessatti, \emph{Sulle superfici di Jacobi semplicemente singolari},
Mem. della societ\`{a} italiana delle Scienze (detta dei XL), 21(3),
45-71, (1919).

\bibitem{DS}W.Decker and F.O.Schreyer, \emph{On the uniqueness of the
Horrocks-Mumford bundle}, Math.Ann. 273 (1986), 415-443.

\bibitem{DM}C.Dionisi and M.Maggesi, \emph{Minimal resolution of general stable
vector bundles on $\PP^2$}, Boll.Un.Mat.It., Sez.B, 6 (2003) 1,
151-160.

\bibitem{Ha74}R.Hartshorne, \emph{Varieties of small codimension in projective
space}, BAMS 80 (1974).

\bibitem{HL}A.Hirschowitz and Y.Laszlo, \emph{Fibrés génériques sur le plan
projectif}, Math. Ann. 297 (1993), 85-102.

\bibitem{Hor}G.Horrocks, \emph{Vector bundles on the punctured spectrum of a
local ring}, Proc.Lond.Math.Soc., III. Ser. 14 (1964), 689-713.

\bibitem{HM}G.Horrocks and D.Mumford, \emph{A rank 2 vector bundle
on $\PP^4$ with 15.000 symmetries}, Topology 12 (1973), 63-81.

\bibitem{Hu}K.Hulek, \emph{The Horrocks-Mumford bundle}, in \emph{Vector bundles in
algebraic geometry (Durham 1993)}, London Math. Soc. Lecture Note
Ser., 208, Cambridge Univ. Press, Cambridge (1995), 139-177.

\bibitem{Ma}N.Manolache, \emph{Syzygies of abelian surfaces embedded in $\PP^4(\C)$},
J. reine angew. Math. 384 (1988), 180-191.

\bibitem{OSS}C.Okonek, M.Schneider and H.Spindler, \emph{Vector bundles on complex
projective spaces}, Progress in Mathematics, no.3, Birkh\"{a}user,
Boston-Basel-Stuttgart, 1980.

\bibitem{OV}G.Ottaviani and J.Vallès, \emph{Moduli of vector bundles and group
action, Notes for Summer School in Algebraic Geometry}, Wykno, 2001.

\bibitem{Va}J.Vallès, \emph{Nombre maximal d'hyperplans instables pour un fibré de
Steiner}, Math. Z. 233 (2000), no. 3, 507-514.

\bibitem{Mac2}D.R. Grayson and M.E. Stillman, \emph{Macaulay 2, a software system for research
in algebraic geometry}, available at http://www.math.uiuc.edu/Macaulay2/.

\end{thebibliography}
\end{document}